\documentclass[11pt,reqno]{amsart}
\usepackage[paperwidth=170mm,paperheight=240mm,left=23mm,right=23mm,top=25mm,bottom=23mm]{geometry}
\usepackage{mathtools,amssymb,amsthm,mathrsfs,microtype}
\usepackage{tikz}
\usepackage[hidelinks]{hyperref}
\usepackage[nameinlink,capitalize,noabbrev]{cleveref}

\newtheorem{theorem}{Theorem}[section]
\newtheorem{lemma}[theorem]{Lemma}
\newtheorem{proposition}[theorem]{Proposition}
\newcommand{\C}{\mathbb C}
\newcommand{\R}{\mathbb R}
\newcommand{\Z}{\mathbb Z}
\newcommand{\F}{\mathbb F}
\newcommand{\1}{\mathbf 1}
\newcommand{\cE}{\mathcal E}
\newcommand{\cR}{\mathcal R}
\newcommand{\tr}{\operatorname{tr}}
\newcommand{\im}{\operatorname{Im}}
\newcommand{\dd}{\,\mathrm d}

\title[Obstructed Harvey--Lawson cones]{Obstructed Harvey--Lawson cones and logarithmic special Lagrangian singularities}
\author[W.-C. Chen]{Wei-Chun Chen}
\author[C.-K. Lien]{Chun-Kai Lien}
\thanks{Department of Mathematics, National Taiwan University,
  Taipei 10617, Taiwan (both authors).
  \texttt{b14201028@ntu.edu.tw} (Chen), \texttt{d13221002@ntu.edu.tw} (Lien).}
\hypersetup{pdfauthor={Wei-Chun Chen and Chun-Kai Lien}}
\date{}
\begin{document}

\begin{abstract}
We construct two exact embedded special Lagrangian germs with the
same unique multiplicity-one tangent cone, the Harvey--Lawson cone
in \(\C^9\). A quadratic obstruction and a corrected cubic obstruction
yield the distinct sharp \(C^1\) convergence rates
\(\lvert\log r\rvert^{-1}\) and \(\lvert\log r\rvert^{-1/2}\),
respectively. These rates distinguish the germs under unitary
transformations and dilations. We also prove Jacobi nonintegrability
of the Harvey--Lawson cones in \(\C^8\) and \(\C^9\), and realize
the sharp rate \(\lvert\log r\rvert^{-1/2}\) at the eight-dimensional cone.
\end{abstract}

\maketitle

\section{Introduction}

Let \(C_{HL}^m\subset\C^m\) be the Harvey--Lawson special
Lagrangian cone \cite{HarveyLawson}. Up to a fixed common phase
rotation, the cone and its link are parametrized by
\begin{gather*}
 \Sigma_m=\{y\in(\R/2\pi\Z)^m:\textstyle\sum_jy_j=0\bmod 2\pi\},\\
 X_0(y)=m^{-1/2}(e^{iy_1},\dots,e^{iy_m}),\qquad
 C_{HL}^m=\{rX_0(y):r\geq 0,\ y\in\Sigma_m\}.
\end{gather*}
Here and below we identify \(\Sigma_m\) with the link by \(X_0\), and
fix the special Lagrangian phase to be that of \(C_{HL}^m\).

The linearized special Legendrian equation in the graph coordinates
of Section~\ref{sec:cone-chart} is \(L_m f=0\), where
\[
 L_m:=\tfrac12(\Delta_{\Sigma_m}-2m)
\]
and \(\Delta_{\Sigma_m}\) is the nonnegative link Laplacian.
In Collins and Li's kernel calculation for \(m\geq3\), the
\(SU(m)\)-rotation modes exhaust the kernel except in dimensions
eight and nine, where there are \(70\) and \(168\) additional real
Jacobi modes, respectively \cite[Section~2.1]{CollinsLi}.
Joyce conjectured that these two cones are Jacobi nonintegrable
\cite[Section~3.2]{JoyceII}.

\subsection{Non-integrability of the cones}

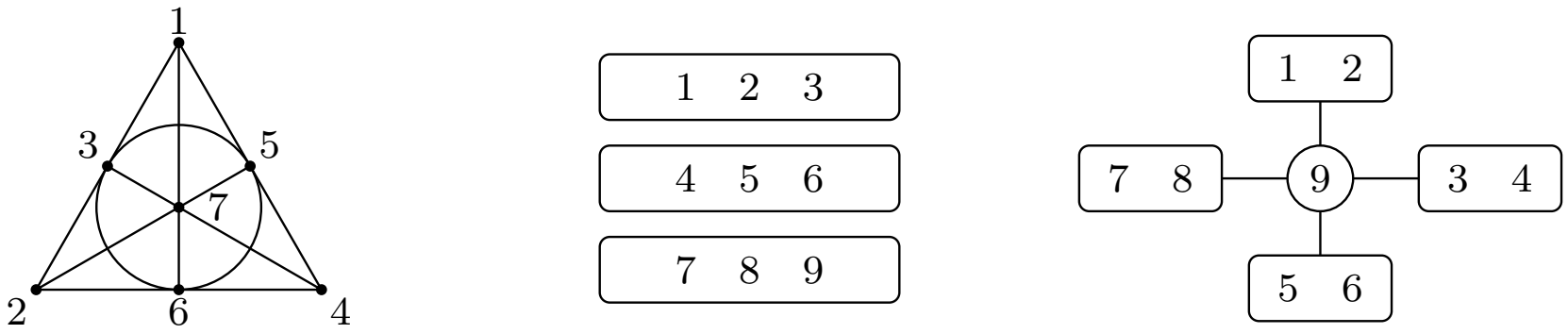
\begin{figure}[!ht]
\centering
\begin{tikzpicture}[line width=.45pt,
    every node/.style={font=\scriptsize,inner sep=1.5pt},
    group/.style={draw,rounded corners=2pt,fill=white,minimum height=4.6mm}]
  \begin{scope}[yshift=-.78cm]
    \coordinate (p1) at (0,{sqrt(3)});
    \coordinate (p2) at (-1,0);
    \coordinate (p4) at (1,0);
    \coordinate (p3) at (-.5,{sqrt(3)/2});
    \coordinate (p5) at (.5,{sqrt(3)/2});
    \coordinate (p6) at (0,0);
    \coordinate (p7) at (0,{sqrt(3)/3});
    \draw (p1)--(p2)--(p4)--cycle
          (p1)--(p6) (p2)--(p5) (p4)--(p3);
    \draw (p7) circle ({sqrt(3)/3});
    \foreach \i/\pos in {1/above,2/below left,3/above left,
        4/below right,5/above right,6/below}
      \fill (p\i) circle[radius=1.1pt] node[\pos] {\(\i\)};
    \fill (p7) circle[radius=1.1pt] node[right=4pt] {\(7\)};
  \end{scope}
  \begin{scope}[xshift=4cm]
    \foreach \height/\triple in {.64/{1\quad 2\quad 3},
        0/{4\quad 5\quad 6},-.64/{7\quad 8\quad 9}}
      \node[group,minimum width=21mm] at (0,\height) {\(\triple\)};
  \end{scope}
  \begin{scope}[xshift=8cm]
    \foreach \x/\y/\pair in {0/.77/{1\quad 2},1.19/0/{3\quad 4},
        0/-.77/{5\quad 6},-1.19/0/{7\quad 8}} {
      \draw (0,0)--(\x,\y);
      \node[group,minimum width=10mm] at (\x,\y) {\(\pair\)};
    }
    \node[draw,circle,fill=white,minimum size=4.6mm] at (0,0) {\(9\)};
  \end{scope}
  \node at (0,-1.3) {Fano configuration};
  \node at (4,-1.3) {Triadic configuration};
  \node at (8,-1.3) {Sunflower configuration};
\end{tikzpicture}
\caption{The configurations for \(u_F,u_T,u_S\), from left to right.
The Fano sets are the complements in \(\{1,\ldots,7\}\) of the seven
lines shown. The triadic boxes are disjoint triplets. Each sunflower
triplet consists of a boxed pair and the common coordinate \(9\).}
\label{fig:configurations}
\end{figure}

Figure~\ref{fig:configurations} shows the three configurations used
below. In dimension eight, the quadratic term vanishes on the excess
Jacobi space, and the Fano configuration gives a cubic obstruction.
In dimension nine, the triadic configuration gives a quadratic
obstruction and the sunflower configuration gives a corrected cubic
obstruction.

McLean proved that deformations of compact smooth special Lagrangians
are unobstructed, with infinitesimal deformations parametrized by
harmonic one-forms \cite[Theorem~3.6]{McLean}. For a special Lagrangian
cone, nonlinear terms in the link equation can obstruct the integration
of a Jacobi field. In Joyce's sense \cite[Definition~6.7]{JoyceI}, \(C_{HL}^m\) is
\emph{Jacobi integrable} if every \(u\in\ker L_m\) is the initial
velocity, in the graph coordinates, of a smooth curve of special
Legendrian links through \(\Sigma_m\).

\begin{theorem}\label{thm:nonint}
The cones \(C_{HL}^8\) and \(C_{HL}^9\) are Jacobi nonintegrable.
\end{theorem}

\subsection{Logarithmic singularities}

The obstructions above give rise to three germs with sharp logarithmic
convergence rates.

\begin{samepage}
\begin{theorem}\label{thm:slow}
There exist exact embedded special Lagrangian germs \(L_T\) and
\(L_S\) sharing the tangent cone \(C_{HL}^9\), and a third such germ
\(L_F\) with tangent cone \(C_{HL}^8\). Each germ is smooth away
from the origin and has a unique tangent cone of multiplicity one.

For each \(m\)-dimensional germ \(L_\ell\), \(\ell\in\{F,T,S\}\),
there are \(c,C>0\) such that for all sufficiently small \(r\), the
rescaling \(r^{-1}L_\ell\) has a normal graph \(x\mapsto x+\nu_{\ell,r}(x)\)
over \(C_{HL}^m\cap\{x:1/2<|x|<2\}\) satisfying
\begin{equation}\label{eq:rate}
 c|\log r|^{-q_\ell}\leq \|\nu_{\ell,r}\|_{C^1}\leq C|\log r|^{-q_\ell},
 \qquad q_F=q_S=\tfrac12,\quad q_T=1.
\end{equation}
In particular, \(L_T\) and \(L_S\) are
inequivalent under unitary transformations and dilations.
\end{theorem}
\end{samepage}
\pagebreak[3]

Adams and Simon constructed minimal submanifolds with logarithmic
convergence \cite{AdamsSimon}, and Joyce anticipated this behavior for
special Lagrangians \cite[Section~6.3]{JoyceI}. Dimler and Gaia construct
exact special Lagrangian perturbations with power decay at any regular
special Lagrangian cone, including \(C_{HL}^9\)
\cite[Theorem~3]{DimlerGaia}. The lower bound in \eqref{eq:rate} rules out \(O(r^\alpha)\) decay of
the rescaled normal graphs for every \(\alpha>0\). These germs thus
fall outside the positive-power setting of Joyce's theory of conical
singularities \cite[Definition~3.7]{JoyceV}.

Let \(\xi(\tau)\) denote the amplitude of the chosen Jacobi mode, with
\(\tau=-\log r\). In each symmetry class, the first nonzero reduced term
\(\kappa\xi^k\) on the critical line gives the leading radial balance,
with \(\kappa>0\) for the branches considered here.
\begin{equation}\label{eq:intro-balance}
 \frac{m+2}{2}\xi'(\tau)+\kappa\xi(\tau)^k=0.
\end{equation}
Its positive solutions decay like \(\lvert\log r\rvert^{-1/(k-1)}\).
The quadratic term for \(L_T\) and the corrected cubic term for \(L_S\)
give \(k=2\) and \(k=3\), respectively. The sharp logarithmic rate
can therefore vary even when the unique tangent cone is fixed.

Simon proved uniqueness for stationary integral varifolds with a
smooth multiplicity-one tangent cone \cite[Corollary, p.~564]{Simon}.
Related results include power decay for integrable minimal cones
\cite{AllardAlmgren} and logarithmic bounds for almost area-minimizing
currents \cite{EngelsteinSpolaorVelichkov}.

The integrability condition of Collins and Li
\cite[Definition~2.2 and Remark~2.6]{CollinsLi}, used in their uniqueness
theorem for cylindrical tangent cones \cite[Theorem~1.1]{CollinsLi},
implies Joyce's condition. Theorem~\ref{thm:nonint} therefore also
rules out their integrability condition for these two cones. Yang Li proved tangent-cone
uniqueness for cones in \(\C^3\) with immersed special Legendrian
links having transverse double points, assuming Jacobi integrability
and a condition on the intersection graph \cite[Theorem~1.2]{YangLi}.

\section*{AI Usage Statement}
Through iterative interactions guided by the authors, 
ChatGPT 5.6 Sol helped identify the symmetries that isolate
the one-dimensional Jacobi kernels in Section~\ref{sec:profiles}
and develop the weighted estimates for the linearized operator
and the nonlinear remainder in Section~\ref{sec:analysis}. The authors independently verified the arguments and take full responsibility for the correctness and content of the manuscript.

\section{The Special Lagrangian Phase Equation}\label{sec:cone-chart}

In the generating-function coordinates below, the phase equation for
deformations of the cone reduces to the special Legendrian equation on
the link when the scaled potential is independent of radius.

\subsection{Exact Lagrangian graphs near the cone}

Let \(\1=(1,\ldots,1)\in\R^m\) and \(\mathcal H=\1^\perp\).
Write \(\Pi:\R^m\to\mathcal H\) for the orthogonal projection. Angular derivatives
\(\nabla f,D^2f\) below use the Euclidean metric on \(\mathcal H\).
The link metric is \(1/m\) times this metric. Set \(r=|z|\), \(R=r^2\),
and \(a=(a_1,\ldots,a_m)\) with \(a_j=|z_j|^2/R\), so \(\sum_j a_j=1\).
Write \(z_j=\sqrt{Ra_j}\,e^{i(y_j+\varphi)}\) with \(y\in\Sigma_m\),
where \(m\varphi\) is the argument of \(z_1\cdots z_m\) in the branch
near zero on a neighborhood of the punctured cone.
The standard Liouville form is
\[
 \lambda=\tfrac R2\sum_j a_j\dd(y_j+\varphi)
        =\tfrac R2\bigl(\dd\varphi+\langle\Pi a,\dd y\rangle\bigr).
\]
For a potential \(G=G(R,y)\), set \(R\Pi a=\nabla_y G\) and
\(\varphi=-G_R\). The Liouville form restricts to
\[
 \lambda=\tfrac12(-R\dd G_R+\dd_yG)=\tfrac12\dd(G-RG_R),
\]
so the graph is exact Lagrangian.
Write \(G=Rf(\tau,y)\) with \(\tau=-\log r\). Then
\(G_R=f-\tfrac12 f_\tau\) and \(\Pi a=\nabla_y f\).
Writing \(r_0>0\) for the radius of the graph neighborhood, the graph
is parametrized by
\begin{equation}\label{eq:conegraph}
 \begin{gathered}
 \bigl\{z_f(r,y):0<r<r_0,\ y\in\Sigma_m\bigr\},\\
 z_f(r,y)=r\,e^{i\varphi}\bigl(\sqrt{a_j}\,e^{iy_j}\bigr)_{j=1}^m,
 \qquad a_j=\frac1m+(\nabla f)_j,\qquad \varphi=-f+\tfrac12 f_\tau .
 \end{gathered}
\end{equation}
For small \(C^1\) jets of \(f\), each \(a_j>0\) and \(m\varphi\) lies in
the chosen branch. The inverse parameters \(r=|z|\) and
\(y_j=\arg z_j-\varphi\pmod{2\pi}\) are smooth, so \eqref{eq:conegraph}
gives an embedded graph.
When \(f\equiv0\), we have \(a_j=1/m\) and \(\varphi=0\), so
\(z_0(r,y)=rX_0(y)\) and the graph is the Harvey--Lawson cone
\(C_{HL}^m\) restricted to \(0<|z|<r_0\).
For static \(f\), the link is a Legendrian graph with \(\varphi=-f\).
The angular projection of a sufficiently \(C^1\)-close Legendrian
link is a local diffeomorphism of degree one. Its contact equation
then gives the graph representation with \(f=-\varphi\).

For the deformation theory of special Legendrian links, see
\cite[Lemma~2.1]{CollinsLi} and \cite{Moriyama}.

\subsection{The special Lagrangian phase}

To compute the Lagrangian angle, factor \(z_j\) out of the \(j\)-th row
of the complex tangent determinant. For \(v\in\mathcal H\), let \(D_v\)
denote the angular derivative in direction \(v\). The column entries satisfy
\[
 \frac{\partial_\tau z_j}{z_j}=-1+\frac{\partial_\tau a_j}{2a_j}+i\varphi_\tau,
 \qquad
 \frac{D_vz_j}{z_j}=\frac{D_va_j}{2a_j}+i(v_j+D_v\varphi).
\]
Let \(\mathfrak D(f)\) be the determinant of these columns in a fixed
orthonormal angular frame. The lifted Lagrangian phase of the graph
relative to the cone is
\[
 \cE_m(f)=m\varphi+\im\log(\mathfrak D(f)/\mathfrak D(0)),
\]
where the logarithm is taken near \(1\).
For static \(f\), the radial column is \(-\1\), so the determinant
uses only the projections of the angular columns onto \(\mathcal H\).
These projections remove the \(D_v\varphi\) term. Factoring out \(i\) from
each projected column and substituting \(a_j\) from
\eqref{eq:conegraph} gives the operator \(F_m\) in the special
Legendrian equation \(F_m(f)=0\)
\begin{equation}\label{eq:Fexact}
 \begin{aligned}
 F_m(f)&=-mf+\im\log\det(I-iM(f)),\\
 M(f)&=\frac m2\Pi\operatorname{diag}\!\bigl((1+m(\nabla f)_j)^{-1}\bigr)
                  D^2 f\big|_{\mathcal H}.
 \end{aligned}
\end{equation}
The phase is analytic in the second jet near zero, with
\[
 \cE_m(tf)=t\bigl(m\varphi-\varphi_\tau-\tfrac m2\tr D^2f\bigr)+O(t^2),
 \qquad \varphi=-f+\tfrac12 f_\tau.
\]
Using \(L_m=\tfrac12(\Delta_{\Sigma_m}-2m)\) and
\(\Delta_{\Sigma_m}=-m\tr D^2\), we obtain the cylinder linearization
\begin{equation}\label{eq:cylinderlinear}
 \mathcal L_m=L_m+\tfrac12(-\partial_\tau^2+(m+2)\partial_\tau).
\end{equation}
For small jets the cylinder equation is uniformly elliptic.
With \(F_m\) acting on each angular slice, define the radial nonlinear
remainder by
\begin{equation}\label{eq:radial-remainder}
 \cE_m(f)=F_m(f)+\tfrac12(-f_{\tau\tau}+(m+2)f_\tau)+\cR_{\rm rad}(f).
\end{equation}

\subsection{The quadratic and cubic terms}

The Taylor expansion of \eqref{eq:Fexact} is
\[
 F_m(tf)=tL_m f+t^2F_{m,2}(f)+t^3F_{m,3}(f)+O(t^4).
\]
Let \(D_j=D_{\Pi e_j}\), where \(e_j\) is the \(j\)-th standard basis vector,
and write \(f_j=D_j f\) and \(f_{jj}=D_j^2 f\). The coefficients are
\begin{align}\label{eq:phase-expansion}
 F_{m,2}(f)&=\frac{m^2}{2}\sum_j f_j f_{jj},\\
 F_{m,3}(f)&=-\frac{m^3}{2}\sum_j f_j^2 f_{jj}
                      +\frac13\tr\!\left(\bigl(\tfrac m2 D^2f\bigr)^3\right).
                      \notag
\end{align}
For the full phase, divide the radial column by \(\sqrt m\) and
use the target frame beginning with \(\1/\sqrt m\). Multiplication on the
right by \(\operatorname{diag}(-1,-iI)\) gives the normalized matrix
\(I+tK_1+t^2K_2+O(t^3)\) for \(tf\), with
\[
 K_1=\begin{pmatrix}
  -i\varphi_\tau & \sqrt m\,(\nabla\varphi)^{\mathsf T}\\
  -\tfrac{\sqrt m}{2}\nabla f_\tau & -i\tfrac m2 D^2f
 \end{pmatrix},
 \qquad \varphi=-f+\tfrac12 f_\tau.
\]
The diagonal blocks are purely imaginary and the off-diagonal blocks
are real, so \(\tr K_1^2\) is real. In the normalized radial column the
imaginary part is linear in \(\varphi_\tau\). In a normalized angular
column, all radial dependence is in the linear first-entry term
\(\sqrt m D_v\varphi\). Hence \(\im\tr K_2\) equals its static value.
The quadratic coefficient \(\im\tr(K_2-\tfrac12 K_1^2)\) is therefore
\(F_{m,2}\), since \(m\varphi\) is linear.

Thus \(\cR_{\rm rad}\) consists of terms of degree at least three with
at least one radial derivative. Write \(j^2f\) for the full jet and
\(j_{\rm rad}^2f=(f_\tau,f_{\tau\tau},\nabla f_\tau)\). For small jets,
\[
\begin{aligned}
 |\cR_{\rm rad}(f)|&\le C|j^2 f|^2|j_{\rm rad}^2f|,\\
 |D\cR_{\rm rad}(f)[g]|&\le C\bigl(
   |j^2 f||j_{\rm rad}^2 f||j^2 g|+|j^2 f|^2|j_{\rm rad}^2g|\bigr),\\
 |D^2\cR_{\rm rad}(f)[g,h]|&\le C|j^2 f||j^2 g||j^2 h|.
\end{aligned}
\]
These bounds follow by differentiating the analytic jet expansion and
also hold with local \(C^{0,\alpha}\) norms in place of absolute values.

The symmetric polarization is \(F_{m,2}(f,g)=\tfrac{m^2}{4}\sum_j
D_j(f_jg_j)\). Integration by parts gives the symmetric identity
for the \(L^2(\Sigma_m)\) pairing
\[
 \langle h,F_{m,2}(f,g)\rangle
 =-\frac{m^2}{4}\int_{\Sigma_m}\sum_j h_j f_j g_j\,\dd\mathrm{vol}_{\Sigma_m}.
\]

\section{Non-Integrability of the Harvey--Lawson Cones}

\subsection{The critical modes}\label{sec:critical-modes}

For \(v\in\Pi\Z^m\), write \(\chi_v(y)=e^{iv\cdot\tilde y}\), where
\(\tilde y\in\mathcal H\) is a zero-sum lift of \(y\in\Sigma_m\).\footnote{Choose
real representatives of the angles \(y_j\). Their sum is \(2\pi k\) for
some \(k\in\Z\), so subtracting \(2\pi k\) from one representative makes
the sum zero. Two such lifts differ by \(2\pi n\) for
\(n\in\Z^m\cap\mathcal H\). If \(v=\Pi n'\) with \(n'\in\Z^m\), then
\(v\cdot n=n'\cdot n\in\Z\), so the exponential is independent of the
lift.}
These characters satisfy \(\Delta_{\Sigma_m}\chi_v=m|v|^2\chi_v\)
\cite[Section~3.2]{JoyceII}. For further spectral calculations see
\cite[Section~6.3.4]{Marshall} and \cite{Haskins}.

For \(A\subset\{1,\ldots,m\}\), define its indicator vector \(\1_A\in\R^m\) by
\[
 (\1_A)_j=
 \begin{cases}
  1, & j\in A,\\
  0, & j\notin A.
 \end{cases}
\]
Set
\[
 v_A=\Pi\1_A,\qquad \theta_A=v_A\cdot\tilde y=\sum_{j\in A}\tilde y_j\pmod{2\pi}.
\]
The rotation space \(\mathcal R_m\) is spanned over \(\R\) by the
real and imaginary parts of the characters with frequencies \(e_i-e_j\),
as follows from the \(SU(m)\) moment maps. The identity
\(m|v_A|^2=|A|(m-|A|)\) gives \(\Delta_{\Sigma_m}\chi_{v_A}=2m\chi_{v_A}\)
for four-element subsets in dimension eight and triplets in dimension nine.
Let \(E_8\) and \(E_9\) be the respective real spans of the cosines
and sines of these \(\theta_A\). Complementation pairs the
four-element subsets in dimension eight, and in dimension nine
each triplet contributes two real modes. These spaces are
orthogonal to the rotation space and have dimensions \(70\) and
\(168\), exhausting the excess multiplicities in
\cite[Section~2.1]{CollinsLi}. Thus in dimensions eight and nine we
have \(\ker L_m=\mathcal R_m\oplus E_m\).

\subsection{Quadratic terms and the triadic configuration}

Extend \(F_{m,2}(f,g)\) complex bilinearly. On Fourier characters
this gives
\begin{equation}\label{eq:quadratic-mixed}
 2F_{m,2}(\chi_v,\chi_w)
 =-\frac{im^2}{2}\sum_j(v_jw_j^2+w_jv_j^2)\,\chi_{v+w}.
\end{equation}
For subsets \(A,B\), the frequency inner product is
\(v_A\cdot v_B=|A\cap B|-|A||B|/m\). The subset formulas are
\[
 \begin{aligned}
 2F_{m,2}(\chi_{v_A},\chi_{v_B})
   &=-im(m-|A|-|B|)(v_A\cdot v_B)\,\chi_{v_A+v_B},\\
 2F_{m,2}(\chi_{v_A},\chi_{-v_B})
   &=-im(|A|-|B|)(v_A\cdot v_B)\,\chi_{v_A-v_B}.
 \end{aligned}
\]

The triadic configuration partitions the nine coordinates into
\(A=\{1,2,3\}\), \(B=\{4,5,6\}\), and \(C=\{7,8,9\}\)
as in Figure~\ref{fig:configurations}. Set
\begin{equation}\label{eq:uT}
 u_T=-\sin\theta_A-\sin\theta_B-\sin\theta_C.
\end{equation}
The frequency vectors satisfy \(v_A+v_B+v_C=0\). Each has squared
length two, and their pairwise inner products are \(-1\).
Among the sums of two frequencies of \(u_T\), only
\(\pm v_A,\pm v_B,\pm v_C\) have squared length two. The only
unordered pair summing to \(-v_C\) is \(\{v_A,v_B\}\), and the first
subset formula gives
\[
 2F_{9,2}\bigl(\tfrac i2\chi_{v_A},\tfrac i2\chi_{v_B}\bigr)
 =27i(i/2)^2\chi_{-v_C}=-\tfrac{27i}{4}\chi_{-v_C}.
\]
Cyclic permutation and complex conjugation give the orthogonal projection
\begin{equation}\label{eq:triadic-coeff}
 \operatorname{proj}_{\ker L_9}F_{9,2}(u_T)
 =-\tfrac{27}{2}\bigl(\sin\theta_A+\sin\theta_B+\sin\theta_C\bigr)
 =\tfrac{27}{2}u_T.
\end{equation}
For \(m=8\) and \(|A|=|B|=4\), the factors \(m-|A|-|B|\) and
\(|A|-|B|\) both vanish. The conjugate interactions vanish as well. Thus
\begin{equation}\label{eq:F2zero}
 F_{8,2}(v,w)=0\qquad(v,w\in E_8).
\end{equation}

\subsection{The Fano configuration}

With the quadratic term zero on \(E_8\), we use the Fano
configuration\footnote{The sets in \(\mathscr F\) are the four-element
complements in \(\{1,\ldots,7\}\) of the lines in
Figure~\ref{fig:configurations}.}
\[
 \mathscr F=\left\{
 \begin{gathered}
  \{1,2,4,7\},\ \{1,2,5,6\},\ \{1,3,4,6\},\ \{1,3,5,7\},\\
  \{2,3,4,5\},\ \{2,3,6,7\},\ \{4,5,6,7\}
 \end{gathered}
 \right\}.
\]
In dimension eight, \(v_A\cdot v_B=|A\cap B|-2\)
for four-element subsets. Any two distinct sets in \(\mathscr F\)
meet in two elements, so the vectors indexed by \(\mathscr F\)
form an orthogonal basis of \(\mathcal H\). Their entries have the same
magnitudes, so a single coefficient calculation gives the cubic
projection onto their cosine sum.

Let \(v_i\in\Pi\Z^m\) be mutually orthogonal frequencies of squared
length two. Write \(\theta_i=v_i\cdot\tilde y\bmod 2\pi\) and
\(V=\sum_i\cos\theta_i\). The cubic term is even in \(y\) and has the
expansion
\begin{equation}\label{eq:orthogonal-cubic}
 F_{m,3}(V)
 =\sum_i\left(
 -\frac{m^3}{4}+\frac{m^3}{8}\sum_j (v_i)_j^4
 +\frac{m^3}{4}\sum_{\ell\ne i}\sum_j (v_i)_j^2(v_\ell)_j^2
 \right)\cos\theta_i+\cdots,
\end{equation}
where the omitted modes have frequencies outside \(\{\pm v_i\}_i\).
Since \(D^2V\,v_i=-2\cos\theta_i\,v_i\), the trace term in
\eqref{eq:phase-expansion} is
\(-\tfrac{m^3}{3}\sum_i\cos^3\theta_i=-\tfrac{m^3}{4}V+\cdots\).
By independence of the frequencies, the other terms at \(\pm v_i\)
follow from
\[
 \begin{aligned}
 \cos\theta_i\sin^2\theta_i&=\tfrac14(\cos\theta_i-\cos3\theta_i),\\
 \cos\theta_i\sin^2\theta_\ell&=\tfrac12\cos\theta_i(1-\cos2\theta_\ell)
 \qquad (\ell\ne i).
 \end{aligned}
\]

For the Fano configuration set
\begin{equation}\label{eq:uF}
 u_F=7^{-1/2}\sum_{A\in\mathscr F}\cos\theta_A.
\end{equation}
The vectors \(v_A=\1_A-\frac12\1\) have squared length two and
entries \(\pm1/2\), so
\[
 \sum_j(v_A)_j^4=\frac12,\qquad
 \sum_j(v_A)_j^2(v_B)_j^2=\frac12\quad(A\ne B).
\]
Taking \(V=\sqrt7\,u_F\) in \eqref{eq:orthogonal-cubic} gives
\[
 F_{8,3}(V)=(-128+32+384)V+\cdots=288V+\cdots.
\]

Cubic homogeneity and the normalization in \eqref{eq:uF} give
\begin{equation}\label{eq:kappa}
 \operatorname{proj}_{\R u_F}F_{8,3}(u_F)=\tfrac{288}{7}u_F.
\end{equation}

\begin{proof}[Proof of Theorem~\ref{thm:nonint}]
A smooth curve of special Legendrian links through \(\Sigma_m\) is
represented by \(f(t)\) with \(f(0)=0\) and \(F_m(f(t))=0\).
For \(u_T\), a solution curve
\(f(t)=tu_T+t^2g+O(t^3)\) would satisfy
\(L_9g+F_{9,2}(u_T)=0\). Since \(L_9\) is self-adjoint,
projection onto its kernel contradicts \eqref{eq:triadic-coeff}.

For \(u_F\), a curve \(f(t)=tu_F+t^2g+t^3h+O(t^4)\) would satisfy
\(L_8g=-F_{8,2}(u_F)=0\) by \eqref{eq:F2zero}. At cubic order it would
satisfy
\[
 L_8h+2F_{8,2}(u_F,g)+F_{8,3}(u_F)=0.
\]
Pairing with \(u_F\) eliminates the first term by self-adjointness.
The middle term vanishes by symmetry of the quadratic pairing,
\[
 \langle u_F,F_{8,2}(u_F,g)\rangle=\langle g,F_{8,2}(u_F)\rangle=0.
\]
The remaining term is \(\tfrac{288}{7}\|u_F\|_{L^2}^2\) by
\eqref{eq:kappa}, a contradiction.
Each cone therefore has a Jacobi direction that fails to integrate.
\end{proof}

\section{Reduced Equations and Asymptotic Profiles}\label{sec:profiles}

A radially constant scaled potential \(f\) deforms the cone itself.
To construct the singular germs, we let the Jacobi amplitude depend on
\(\tau\). The resulting first-order radial term in
\eqref{eq:cylinderlinear} balances the static obstruction at leading
order. For each direction we restrict the phase equation to a symmetry
class whose Jacobi kernel is spanned by that direction.

\subsection{Symmetry classes and critical lines}

In dimension eight, the Fano configuration also determines the
translations that select its seven frequency pairs among the Jacobi
modes. Write \(x=(y_1,\ldots,y_7)\), so
\(y_8=-\sum_a x_a\) modulo \(2\pi\). Label the Fano points by the
nonzero elements of \(\F_2^3\). Evaluating linear functionals there
gives a binary code \(S_F\subset\F_2^7\). The three nonzero points in the
kernel of a nonzero functional form a Fano line, so the seven nonzero
words of \(S_F\) are the indicators of the sets in \(\mathscr F\).
Let \(S_F^\perp\) be the orthogonal complement for the standard dot
product over \(\F_2\). Its elements \(h\) are represented
by vectors in \(\{0,1\}^7\) in the translations below.
Let \(\mathcal X_F\) be the space of real angular functions fixed by
translations \(x\mapsto x+\pi h\) for \(h\in S_F^\perp\), permutations
preserving the Fano lines, and the involution
\(f(x)\mapsto-f((\pi/4,\ldots,\pi/4)-x)\).

For the triadic configuration in \eqref{eq:uT}, let
\(\mathcal X_T\) consist of angular functions invariant
under translations of sum zero within each block and cyclic
permutation of the blocks, and odd under \(y\mapsto-y\).

The sunflower configuration in Figure~\ref{fig:configurations} has
four triplets \(A_i=\{2i-1,2i,9\}\) for \(1\leq i\leq4\), all sharing
coordinate \(9\).
To select the cosine direction, choose \(t\in\Sigma_9\) with
\(\theta_{A_i}(t)=\pi\) for every \(i\), for example
\(t=(-\pi/6,\ldots,-\pi/6,4\pi/3)\).
Let \(\mathcal X_S\) be fixed by opposite translations
within each pair \(\{2i-1,2i\}\), permutations of the pairs, and
\(f(y)\mapsto-f(t-y)\). Set
\begin{equation}\label{eq:uS}
 u_S=\tfrac12\sum_{i=1}^4\cos\theta_{A_i}.
\end{equation}

\begin{lemma}\label{lem:symmetry-kernel}
The phase operators preserve the three symmetry classes, with
the angular symmetries imposed on each slice for \(\cE_m\).
Their critical eigenspaces are
\[
\begin{aligned}
 \ker L_8\cap\mathcal X_F&=\R u_F,\\
 \ker L_9\cap\mathcal X_T&=\R u_T,\\
 \ker L_9\cap\mathcal X_S&=\R u_S.
\end{aligned}
\]
\end{lemma}

\begin{proof}
Both phase operators are equivariant under angular translations
and coordinate permutations. Complex conjugation acts on both
potentials and phase forcing by \(f(y)\mapsto-f(-y)\), since it
negates \(\varphi\) and conjugates the determinant ratio at the
reflected point. Composing with translations gives the two
shifted involutions.

For \(n\in\Z^7\), translation by \(\pi h\) multiplies \(e^{in\cdot x}\)
by \((-1)^{n\cdot h}\), so
\[
 e^{in\cdot(x+\pi h)}=e^{in\cdot x}\ \text{for every }h\in S_F^\perp
 \quad\Longleftrightarrow\quad n\bmod 2\in S_F.
\]
A rotation phase \(y_i-y_j\) with \(i,j\leq7\) has two odd coefficients,
and \(y_i-y_8=x_i+\sum_a x_a\) has six. Both weights are excluded.
Each excess pair has one four-element representative \(A\) avoiding
coordinate \(8\), by complementation. Its parity is the indicator of
\(A\), so the seven nonzero words of \(S_F\) select precisely the pairs
from \(\mathscr F\). The action of \(GL(3,\F_2)\) is transitive on
the Fano lines and hence on these sets, leaving the cosine and sine sums. The shifted involution sends
\(\theta_A\) to \(\pi-\theta_A\). With the sign on \(f\), it fixes the
cosine sum and negates the sine sum, leaving exactly \(\R u_F\).

Blockwise translations in \(\mathcal X_T\) select
the three triplet pairs at the critical eigenvalue. Cyclic
permutation and oddness fix the line of \(u_T\).

In \(\mathcal X_S\), a representative constant on each pair
can only be critical at one of the \(\pm v_{A_i}\). Pair permutations
leave the cosine and sine sums. At \(\theta_{A_i}(t)=\pi\) the signed
involution fixes the cosine sum and negates the sine sum.
\end{proof}

In each symmetry class, let \(P\) be the \(L^2(\Sigma_m)\)-orthogonal
projection onto its critical line and set \(Q=I-P\) on that class.
The restriction \(L_m|_Q\) to the range of \(Q\) has a spectral gap
at zero.

\subsection{The sunflower configuration and its quadratic correction}

For Fano, \(F_{8,2}(u_F)=0\). For sunflower, the quadratic term is
nonzero and lies in the complement of the critical space. Its
cancellation contributes to the cubic coefficient.
The vectors \(v_{A_i}\) have squared length two. For \(i\ne\ell\),
the triplets meet at one coordinate, so \(v_{A_i}\cdot v_{A_\ell}=0\).
The subset formulas following \eqref{eq:quadratic-mixed} then give
zero mixed quadratic interactions.
Since \(\sum_j(v_{A_i})_j^3=2/3\),
\(F_{9,2}(g(\theta_{A_i}))=27g'g''\). In particular,
\[
 F_{9,2}(u_S)=\tfrac{27}{8}\sum_i\sin 2\theta_{A_i} .
\]
Since \(L_9\sin 2\theta_{A_i}=27\sin 2\theta_{A_i}\), the forcing is orthogonal
to the full kernel of \(L_9\). To cancel the quadratic term in
\(F_9(\xi u_S+\xi^2h_S)\), set
\begin{equation}\label{eq:hS}
 h_S=-\tfrac18\sum_i\sin 2\theta_{A_i},\qquad L_9h_S=-F_{9,2}(u_S).
\end{equation}
In particular \(h_S\in\mathcal X_S\) and \(Ph_S=0\).

The direct cubic projection uses \eqref{eq:orthogonal-cubic} with
the fourth moments
\(\sum_j(v_A)_j^4=2/3\) and \(\sum_j(v_A)_j^2(v_B)_j^2=4/9\)
for distinct sunflower triplets \(A,B\).
For the mixed term, symmetry of the quadratic pairing gives
\[
 2\langle u_S,F_{9,2}(u_S,h_S)\rangle
 =-2\langle h_S,L_9h_S\rangle=-\tfrac{27}{8}\|u_S\|_{L^2}^2,
\]
using \(L_9h_S=27h_S\) and \(\|h_S\|_{L^2}^2=\|u_S\|_{L^2}^2/16\).
The two projections are therefore
\[
 P F_{9,3}(u_S)=\tfrac{243}{8}u_S,\qquad
 2P F_{9,2}(u_S,h_S)=-\tfrac{27}{8}u_S.
\]
Their sum is
\begin{equation}\label{eq:sunflower-cubic}
 P\bigl(F_{9,3}(u_S)+2F_{9,2}(u_S,h_S)\bigr)
 =27u_S.
\end{equation}
Any solution curve tangent to \(u_S\) would have quadratic term
\(h_S+g\) with \(g\in\ker L_9\). The quadratic pairing gives
\(\langle u_S,F_{9,2}(u_S,g)\rangle=\langle g,F_{9,2}(u_S)\rangle=0\).
Thus the cubic equation paired with \(u_S\) would have left side
\(27\|u_S\|_{L^2}^2\), proving the obstruction among all links.

\subsection{The reduced static equation}

On each critical line set \(\Phi(\xi)=\xi u+\xi^2h\), using
\(h=h_S\) for \(u=u_S\) and \(h=0\) on the other two
lines. Since \(Ph=0\), \(P\Phi(\xi)=\xi u\).
Equations~\eqref{eq:kappa}, \eqref{eq:triadic-coeff}, and
\eqref{eq:sunflower-cubic} give
\begin{equation}\label{eq:static-profile-order}
 \mathcal S(\xi):=F_m(\Phi(\xi))=O_{C^{0,\alpha}}(|\xi|^k),\qquad
 P\mathcal S(\xi)=\kappa\xi^k u+O(|\xi|^{k+1}).
\end{equation}
Here \((k,\kappa)=(3,288/7)\) for \(u=u_F\), \((2,27/2)\) for \(u=u_T\),
and \((3,27)\) for \(u=u_S\). In the cubic sectors, the quadratic
correction cancels the entire second-order static term.
\begin{equation}\label{eq:full-cancellation}
 \mathcal S''(0)=2\bigl(L_mh+F_{m,2}(u)\bigr)=0
 \qquad (k=3).
\end{equation}
In the quadratic sector \(k=2\), \(P\mathcal S''(0)=2\kappa u\ne0\).
Since \(L_m|_Q\) is invertible and \(Q\mathcal S(\xi)=O(\xi^k)\),
the implicit function theorem gives a correction \(\eta(\xi)=O(\xi^k)\)
with \(P\eta(\xi)=0\) and \(QF_m(\Phi(\xi)+\eta(\xi))=0\).
As \(PL_m=0\), the correction changes \(P\mathcal S(\xi)\) only at
order \(k+1\) and higher. For the cylinder equation, we use \(\Phi\)
as the approximate profile and solve for the complementary correction
\(w\) in Section~\ref{sec:analysis}.

\section{Weighted Inversion and Solvability}\label{sec:analysis}

In each class of Lemma~\ref{lem:symmetry-kernel}, put
\(c=m+2\), \(q=1/(k-1)\), and \(\xi_0=A\tau^{-q}\), with \(A>0\)
determined by \(\kappa A^{k-1}=cq/2\). For large \(T\), we seek a
solution \(f=\Phi(\xi_0+b)+w\) on the half-cylinder \(\tau\ge T\),
with \(Pw=0\) and weights
\[
 \beta=q+\tfrac12,\qquad \gamma=kq=q+1.
\]
The weight \(\gamma\) matches the complementary residual
\(Q\cE_m(\Phi(\xi_0))=O(\tau^{-q-1})\) estimated below. The weight
\(\beta\) gives \(b\) an extra half-power of decay relative to \(\xi_0\).
The leading profile satisfies
\begin{equation}\label{eq:balance}
 \tfrac c2 \xi_0'+\kappa \xi_0^k=0,\qquad
 k\kappa \xi_0^{k-1}=\frac{c\gamma}{2\tau}.
\end{equation}

\subsection{Weighted right inverses}

Fix \(\alpha\in(0,1)\) and write \(I_s=[s,s+1]\) for \(s\geq T\geq2\).
The center norm assigns an extra power of decay to the radial
derivatives used in the estimates for \(\cR_{\rm rad}\).
\begin{align*}
 \|\psi\|_{Y^{\mathrm{sc}}_\sigma}
   &=\sup_{s\geq T}s^\sigma\|\psi\|_{C^{0,\alpha}(I_s)},\\
 \|\psi\|_{Y^{\mathrm{cyl}}_\sigma}
   &=\sup_{s\geq T}s^\sigma\|\psi\|_{C^{0,\alpha}(I_s\times\Sigma_m)},\\
 \|b\|_{B_\beta}
   &=\sup_{s\geq T}\bigl(s^\beta\|b\|_{C^0(I_s)}
             +s^{\beta+1}\|b'\|_{C^{1,\alpha}(I_s)}\bigr),\\
 \|w\|_{Z_\gamma}
   &=\sup_{s\geq T}s^\gamma\|w\|_{C^{2,\alpha}(I_s\times\Sigma_m)}.
\end{align*}
The scalar spaces \(Y^{\mathrm{sc}}_\sigma\) and \(B_\beta\) consist of
real-valued functions on \([T,\infty)\) with finite corresponding norm.
The cylinder spaces \(Y^{\mathrm{cyl}}_\sigma\) and \(Z_\gamma\) consist of
real-valued functions on \([T,\infty)\times\Sigma_m\) in the chosen
symmetry class with finite corresponding norm.
The notation \(QY^{\mathrm{cyl}}_\sigma,QZ_\gamma\) denotes the closed subspaces
with zero \(P\)-projection. Decay orders below are understood in
the corresponding slab H\"older norms.

\begin{lemma}\label{lem:center}
For \(c>0\) and \(p>\beta>0\), the operator
\(\mathcal A_p=\tfrac12(-\partial_\tau^2+c\partial_\tau+cp/\tau)\)
has a right inverse \(J_p:Y^{\mathrm{sc}}_{\beta+1}\to B_\beta\) with norm bounded
independently of sufficiently large \(T\).
\end{lemma}

\begin{proof}
The center operator has the factorization
\[
 2\mathcal A_p=(-\partial_\tau+c)(\partial_\tau+p/\tau)
                    +(p/\tau)\partial_\tau-p/\tau^2 .
\]
For forcing \(\psi\), define \(J_{p,0}\psi:=b\) by
\[
 g(\tau)=2\int_\tau^\infty e^{-c(s-\tau)}\psi(s)\,\dd s,\qquad
 b(\tau)=\tau^{-p}\int_T^\tau s^p g(s)\,\dd s.
\]
These formulas give \(b'+pb/\tau=g\) and \(-g'+cg=2\psi\).
The exponential integral gives \(g=O(\tau^{-\beta-1})\) per
unit forcing norm, and \(p>\beta\) gives
\[
 \tau^{-p}\int_T^\tau s^{p-\beta-1}\,\dd s
 \leq \frac{\tau^{-\beta}}{p-\beta}.
\]
These bounds and the two first-order equations give
\[
 \|J_{p,0}\psi\|_{B_\beta}\leq C\|\psi\|_{Y^{\mathrm{sc}}_{\beta+1}}.
\]
The factorization gives
\[
 \mathcal A_pJ_{p,0}\psi-\psi
 =\frac{p}{2\tau}(J_{p,0}\psi)'-\frac{p}{2\tau^2}J_{p,0}\psi.
\]
The preceding bounds imply \(\|\mathcal A_pJ_{p,0}-I\|\leq CT^{-1}\) on
\(Y^{\mathrm{sc}}_{\beta+1}\), so the Neumann series defines the
uniformly bounded right inverse
\[
 J_p=J_{p,0}(\mathcal A_pJ_{p,0})^{-1}.
\]
\end{proof}

\begin{lemma}\label{lem:complement}
For \(\gamma>0\), \(\mathcal L_m|_Q\) has a right inverse
\(QY^{\mathrm{cyl}}_\gamma\to QZ_\gamma\) with norm bounded independently of \(T\).
\end{lemma}

\begin{proof}
Extend the forcing \(\psi\) to the full cylinder by composing in
time with \(\tau\mapsto\max(T,\tau)\). This preserves the weighted
local H\"older bound with weight \(\omega_T(\tau)=\max(T,\tau)^\gamma\),
including on slabs crossing \(T\). For an eigenvalue
\(\lambda\ne0\) of \(\Delta_{\Sigma_m}-2m\) on the complement,
the characteristic roots are \(r_\pm=(c\pm\sqrt{c^2+4\lambda})/2\).
They are distinct since \(c^2+4\lambda\geq(m-2)^2\).
The Green kernel is
\[
 K_\lambda(t)=\frac{2}{r_+-r_-}
 \begin{cases}
  e^{r_-t}, & \lambda>0,\ t\geq0,\\
  e^{r_+t}, & \lambda>0,\ t<0,\\
  e^{r_+t}-e^{r_-t}, & \lambda<0,\ t<0,\\
  0, & \lambda<0,\ t\geq0 .
 \end{cases}
\]
The kernels are continuous and satisfy \(K_\lambda'(0^+)-K_\lambda'(0^-)=-2\).
Let \(K(t)\) act by \(K_\lambda(t)\) on each complementary eigenspace
of \(\Delta_{\Sigma_m}-2m\) and by zero on the remaining modes.
The spectral gap and the finitely many negative modes give
\[
 \|K(t)\|_{L^2(\Sigma_m)\to L^2(\Sigma_m)}\leq Ce^{-\delta|t|}.
\]
Since \(\omega_T(\tau)/\omega_T(\tau-t)\leq(1+|t|)^\gamma\),
the convolution \(w(\tau)=\int_{\R}K(t)\psi(\tau-t)\,\dd t\)
satisfies
\[
 \|\omega_T w\|_{L^\infty_\tau L^2(\Sigma_m)}
 \leq C\|\psi\|_{Y^{\mathrm{cyl}}_\gamma},
\]
where the forcing norm is taken on \(\tau\geq T\). Testing the
modal equations against finite angular spectral sums and passing
to the limit gives \(\mathcal L_m w=\psi\) in distributions on the
full cylinder.

The weighted \(L^2\) bound and interior \(W^{2,p}\), Sobolev and
Schauder estimates on nested cylinders
\cite[Chapters~6 and~9]{GilbargTrudinger} give
\[
 \|w\|_{Z_\gamma}\leq C\|\psi\|_{Y^{\mathrm{cyl}}_\gamma}.
\]
The constant is uniform in \(T\), since \(\Sigma_m\) is compact and
\(\omega_T\asymp s^\gamma\) on \([s-1,s+2]\) for \(s\geq T\).
\end{proof}

\subsection{Quantitative estimates and solvability}

Use the spaces
\[
 X=B_\beta\oplus QZ_\gamma,\qquad
 Y=Y^{\mathrm{sc}}_{\beta+1}\oplus QY^{\mathrm{cyl}}_\gamma .
\]
For a forcing \(g\), use the slice-wise coordinates
\[
 g_P(\tau)=\frac{\langle g(\tau,\cdot),u\rangle_{L^2(\Sigma_m)}}
                 {\|u\|_{L^2(\Sigma_m)}^2},
 \qquad Pg=g_Pu,\qquad g_Q=Qg.
\]
Write \(\mathscr E(b,w)=(g_P,g_Q)\) for \(g=\cE_m(\Phi(\xi_0+b)+w)\).
Set \(\mathcal D=D\mathscr E(0,0)\) and
\(\mathcal N(v)=\mathscr E(v)-\mathscr E(0,0)-\mathcal Dv\) for
\(v=(b,w)\). In the block notation for \(\mathcal D\), the first subscript
is the output component, so \(\mathcal D_{PQ}w=(\mathcal D(0,w))_P\).
For every fixed \(M\) and all sufficiently large \(T\),
\(\mathcal D\) has a right inverse \(J_{\mathcal D}:Y\to X\) such that
\begin{align}
 \|\mathscr E(0,0)\|_Y&\leq C,\label{eq:approx-error}\\
 \|J_{\mathcal D}\|_{Y\to X}&\leq C,\label{eq:linear-inverse}\\
 \|\mathcal N(v)-\mathcal N(\widetilde v)\|_Y
   &\leq C_M T^{-1/2}\|v-\widetilde v\|_X.
   \label{eq:nonlinear-estimate}
\end{align}
The last estimate holds for \(v,\widetilde v\) in the ball of radius
\(M\). The constants are independent of sufficiently large \(T\),
and the first two bounds are independent of \(M\).

At \(\Phi(\xi_0)\), we have \(\xi_0^{(j)}=O(\tau^{-q-j})\)
for \(0\leq j\leq 2\), so the jet bounds give
\(\cR_{\rm rad}(\Phi(\xi_0))=O(\tau^{-3q-1})\). The leading critical
terms cancel by \eqref{eq:balance}, so \eqref{eq:radial-remainder} and
\eqref{eq:static-profile-order} give
\[
\begin{aligned}
 (\mathscr E(0,0))_P&=O(\tau^{-q-2}+\tau^{-2q-1}),\\
 (\mathscr E(0,0))_Q&=O(\tau^{-q-1}).
\end{aligned}
\]
Since \(q\geq1/2\), these bounds imply \eqref{eq:approx-error}.

In the chosen weights,
\begin{align*}
 \|\mathcal D_{PP}-\mathcal A_\gamma\|_{B_\beta\to Y^{\mathrm{sc}}_{\beta+1}}&\leq CT^{-q},\\
 \|\mathcal D_{QQ}-\mathcal L_m|_Q\|_{QZ_\gamma\to QY^{\mathrm{cyl}}_\gamma}
    &\leq CT^{-q},\\
 \|\mathcal D_{QP}\|_{B_\beta\to QY^{\mathrm{cyl}}_\gamma}&\leq CT^{-1/2}.
\end{align*}
For the center block, the linear radial term on \(b u\) and
\(k\kappa\xi_0^{k-1}b\) give \(\mathcal A_\gamma b\) by \eqref{eq:balance}.
The remaining static term is \(O(s^{-kq-\beta})\) for unit \(b\).
In the first variation of \(\cR_{\rm rad}\), a radial derivative gains
\(s^{-1}\) whether it falls on the profile or on the center variation,
so this term is \(O(s^{-2q-\beta-1})\). The \(QQ\) bound follows from
analyticity and the \(O(s^{-q})\) base jet. For the lower-left block,
the analytic expansion \eqref{eq:static-profile-order} gives \(Q\mathcal S'(\xi_0)=O(\tau^{-1})\),
since \((k-1)q=1\). The radial terms, including derivatives of
\(2\xi_0hb\), are of higher order. The gap
\(\gamma-\beta=\tfrac12\) then gives the \(T^{-1/2}\) bound.
Since \(P\mathcal L_m w=0\) and \(\cR_{\rm rad}\) starts at degree
three, the upper-right block satisfies, on \(I_s\),
\[
 \begin{aligned}
 &\|\mathcal D_{PQ}w-2\xi_0(F_{m,2}(u,w))_P\|_{C^{0,\alpha}(I_s)}\\
 &\qquad\leq Cs^{-2q}\|w\|_{C^{2,\alpha}(I_s\times\Sigma_m)}.
 \end{aligned}
\]
The leading term is controlled for every \(w\in QZ_\gamma\) by
\[
 s^{\beta+1}\|\xi_0\|_{C^{0,\alpha}(I_s)}
             \|w\|_{C^{2,\alpha}(I_s\times\Sigma_m)}
 \le C s^{1/2-q}\|w\|_{Z_\gamma}.
\]
The remainder has weighted size \(O(T^{1/2-2q})\). For triadic, \(q=1\)
gives \(\mathcal D_{PQ}=O(T^{-1/2})\) in operator norm. For Fano, \(q=1/2\)
and the leading term vanishes by the quadratic pairing and
\(F_{8,2}(u_F)=0\), again giving \(O(T^{-1/2})\). For sunflower, \(q=1/2\)
gives a uniform bound \(QZ_\gamma\to Y^{\mathrm{sc}}_{\beta+1}\). The
mixed projection on \(w=\tau^{-\gamma}h_S\) gives
\[
 \mathcal D_{PQ}w=-\tfrac{27}{8}A\tau^{-\beta-1}
                    +O(\tau^{-\beta-3/2}),
\]
since \(q+\gamma=\beta+1\). This \(w\) has uniformly bounded \(Z_\gamma\)
norm, and the weighted value at \(\tau=T\) has a positive lower bound
in absolute value for large \(T\). We therefore retain the sunflower
coupling in the principal operator
\[
 \mathcal T=
 \begin{pmatrix}
  \mathcal A_\gamma & \mathcal T_{PQ}\\
  0 & \mathcal L_m|_Q
 \end{pmatrix},
 \qquad
 \mathcal T_{PQ}=
 \begin{cases}
  \mathcal D_{PQ} & \text{for sunflower},\\
  0 & \text{for Fano and triadic}.
 \end{cases}
\]
With this choice we obtain
\[
 \|\mathcal D-\mathcal T\|_{X\to Y}\leq CT^{-1/2}.
\]
Since \(\gamma>\beta\), Lemma~\ref{lem:center} gives the uniformly
bounded right inverse \(J_\gamma\) of \(\mathcal A_\gamma\). With \(J_Q\)
from Lemma~\ref{lem:complement}, a right inverse of \(\mathcal T\) is
\[
 J_{\mathcal T}(g_P,g_Q)=
 \bigl(J_\gamma(g_P-\mathcal T_{PQ}J_Qg_Q),\,J_Qg_Q\bigr),
\]
with \(\|J_{\mathcal T}\|_{Y\to X}\leq C\). Since
\(\mathcal D J_{\mathcal T}-I=(\mathcal D-\mathcal T)J_{\mathcal T}\),
\begin{equation}\label{eq:linear-error}
 \|\mathcal D J_{\mathcal T}-I\|_{Y\to Y}\leq CT^{-1/2}.
\end{equation}
For large \(T\), the Neumann series defines
\(J_{\mathcal D}=J_{\mathcal T}(\mathcal D J_{\mathcal T})^{-1}\) with the bound
\eqref{eq:linear-inverse}.

For \eqref{eq:nonlinear-estimate}, set \(\xi=\xi_0+b\). Analyticity and
\eqref{eq:full-cancellation} give \(\mathcal S''(\xi)=O(|\xi|^{k-2})\).
At \(w=0\) the chain rule gives
\(D_b^2[F_m(\Phi(\xi))](b_1,b_2)=\mathcal S''(\xi)b_1b_2\).
For unit center variations, \(b_i=O(\tau^{-\beta})\), while
\(\xi=O_M(\tau^{-q})\) and \(w=O_M(\tau^{-\gamma})\).
Analyticity and the H\"older product estimate give
\[
 \begin{aligned}
 &\bigl\|D_b^2[F_m(\Phi(\xi)+w)](b_1,b_2)\bigr\|_{C^{0,\alpha}(I_s\times\Sigma_m)}\\
 &\qquad\leq C_M\bigl(s^{-(k-2)q}+s^{-\gamma}\bigr)s^{-2\beta}
 \leq C_Ms^{-q-2}.
 \end{aligned}
\]
The \(s^{-\gamma}\) term comes from the change from \(w=0\).
The mixed and complementary static variations are
\(O_M(s^{-\beta-\gamma})\) and \(O_M(s^{-2\gamma})\) on unit inputs.
The cubic remainder contributes \(O_M(s^{-q-2\beta})\) by the jet
bounds above. Indeed the first and second variation jets of the
profile map are \(O_M(s^{-\beta})\) and \(O(s^{-2\beta})\), respectively.
The linear radial operator acts on \(2hb_1b_2\), giving
\(O(s^{-2\beta-1})\) by the extra derivative weight in \(B_\beta\).
The center--center static bound has exponent \((k-2)q+2\beta=q+2\).
The other exponents satisfy
\[
 \min\{\beta+\gamma,2\gamma,q+2\beta,2\beta+1\}\geq q+2
\]
because \(q\geq1/2\). Both projections preserve the slab bounds.
Multiplication by \(s^{\beta+1}\) gives \(O_M(s^{-1/2})\) for the
\(P\)-component, and multiplication by \(s^\gamma\) gives
\(O_M(s^{-1})\) for the \(Q\)-component. Consequently
\[
 \sup_{\|v\|_X\leq M}\|D^2\mathscr E(v)\|_{X\times X\to Y}
 \leq C_M T^{-1/2}.
\]
The mean-value theorem and \(D\mathcal N(0)=0\) now give
\eqref{eq:nonlinear-estimate}.

\begin{proposition}\label{prop:sector-solution}
For each of the three sectors and all sufficiently large \(T\), there
are \(b\in B_\beta\) and \(w\in QZ_\gamma\) such that
\(f=\Phi(\xi_0+b)+w\) solves \(\cE_m(f)=0\) on the half-cylinder.
The norms of \(b\) and \(w\) are bounded independently of sufficiently
large \(T\). The solution is smooth on the open half-cylinder.
\end{proposition}

\begin{proof}
By \eqref{eq:approx-error}--\eqref{eq:nonlinear-estimate}, the map
\[
 v\longmapsto-J_{\mathcal D}\bigl(\mathscr E(0,0)+\mathcal N(v)\bigr)
\]
contracts a fixed sufficiently large ball in \(X\) for large \(T\).
The fixed point solves the phase equation, and interior elliptic
bootstrapping gives smoothness on the open half-cylinder
\cite[Chapters~6 and~17]{GilbargTrudinger}.
\end{proof}

\section{Proof of the Geometric Theorems}\label{sec:main-proofs}

\begin{theorem}\label{thm:asymptotics}
The special Lagrangians in Theorem~\ref{thm:slow} may be chosen with
the following rescaled generating functions in the chart of
Section~\ref{sec:cone-chart}, where \(\tau=-\log r\):
\begin{align}\label{eq:mainasymp}
 f_F&=\frac{\sqrt{35}}{24}\tau^{-1/2}u_F+O(\tau^{-1}),\\
 f_T&=\frac{11}{27}\tau^{-1}u_T
                              +O(\tau^{-3/2}),\notag\\
 f_S&=\frac{\sqrt{33}}{18}\tau^{-1/2}u_S
                              +O(\tau^{-1}).\notag
\end{align}
Here \(u_F,u_T,u_S\) are given in
\eqref{eq:uF}, \eqref{eq:uT}, and \eqref{eq:uS}. The errors hold in
\(C^{2,\alpha}\) on unit cylinders for any fixed \(\alpha\in(0,1)\)
in the construction. In the sunflower class, the orthogonal projection
\(Q\) off \(u_S\) gives
\begin{equation}\label{eq:sunflower-Qasymp}
 Q f_S=\frac{11}{108}\tau^{-1}h_S
                          +O(\tau^{-3/2}),
\end{equation}
where \(h_S\) is defined in \eqref{eq:hS}.

For \(\ell\in\{F,T,S\}\), let \(A_\ell\) be the coefficient of
\(\tau^{-q_\ell}u_\ell\) in \eqref{eq:mainasymp}. For \(x=\rho X_0(y)\)
with \(\rho=|x|\), set \(U_\ell(x)=\rho^2u_\ell(y)\).
The section \(\nu_{\ell,r}\) from Theorem~\ref{thm:slow} satisfies
\begin{equation}\label{eq:normal-profile}
 \tau^{q_\ell}\nu_{\ell,r}\longrightarrow -\frac{A_\ell}{2}J\nabla_{C_{HL}^m}U_\ell
 \qquad\text{in }C^1
\end{equation}
on the base annulus, where \(J\) is the ambient complex structure.
\end{theorem}

\begin{proof}[Proof of Theorems~\ref{thm:slow} and \ref{thm:asymptotics}]
Proposition~\ref{prop:sector-solution} and \(\kappa A^{k-1}=cq/2\)
give \eqref{eq:mainasymp}. In the sunflower class,
\(Q f_S=(\xi_0+b)^2h_S+w\), where \(b\) is the scalar correction from
Section~\ref{sec:analysis}. This identity and \(A^2=11/108\) give
\eqref{eq:sunflower-Qasymp}.

The embedded graph \eqref{eq:conegraph} has the cone phase and
primitive \(\tfrac12(G-R\partial_R G)\).
The mass in \(B_r\) is \(O(r^m)\) and the slice mass is \(O(r^{m-1})\), so
the calibrated current extends without boundary at the origin.
The rescalings have mass \(O(\delta^m)\) in \(B_\delta\) and converge
by \eqref{eq:mainasymp} with one sheet in \(C^1\) on compact annuli.
Hence \(C_{HL}^m\) is the unique multiplicity-one tangent cone.
Its nonplanarity makes the origin singular.

Write \(x=\rho X_0(y)\) on a slightly larger annulus and parametrize
the rescaling by \(X_r(x)=r^{-1}z_{f_\ell}(r\rho,y)\).
By \eqref{eq:mainasymp},
\(\tau^{q_\ell}f_\ell(\tau-\log\rho,y)\to A_\ell u_\ell(y)\) in \(C^2\), and
\(\tau^{q_\ell}\partial_\tau f_\ell(\tau-\log\rho,y)\to0\) in \(C^1\).
The first variation of \eqref{eq:conegraph} in the static direction \(u_\ell\) has \(j\)-th
component \(x_j(\tfrac m2 D_j u_\ell-iu_\ell)\). This is the \(j\)-th component
of \(-\tfrac12 J\nabla_{C_{HL}^m}U_\ell\), so
\[
 X_r(x)=x-\tfrac{A_\ell}{2}\tau^{-q_\ell}J\nabla_{C_{HL}^m}U_\ell(x)
           +o_{C^1}(\tau^{-q_\ell}).
\]
The leading field is normal since the cone is Lagrangian. Thus the
tubular footpoint map is the identity up to \(o_{C^1}(\tau^{-q_\ell})\),
and reparametrizing by it gives \eqref{eq:normal-profile} on the base
annulus. The limit is nonzero since \(\partial_\rho U_\ell
=2\rho u_\ell\not\equiv0\) there. Taking \(C^1\) norms in
\eqref{eq:normal-profile} shows that \(\tau^{q_\ell}\|\nu_{\ell,r}\|_{C^1}\)
tends to a finite positive limit, giving both bounds in \eqref{eq:rate}.

An identifying unitary map would preserve the common tangent cone
and the normal graph norms over it. A fixed dilation shifts \(\tau\)
by a constant, so the distinct decay exponents \(1\) and \(1/2\)
rule out equivalence of the two nine-dimensional germs.
\end{proof}


\begin{thebibliography}{99}

\bibitem{AdamsSimon}
D.~Adams and L.~Simon,
\emph{Rates of asymptotic convergence near isolated singularities of
geometric extrema},
Indiana Univ. Math. J. \textbf{37} (1988), no.~2, 225--254.

\bibitem{AllardAlmgren}
W.~K. Allard and F.~J. Almgren, Jr.,
\emph{On the radial behavior of minimal surfaces and the uniqueness
of their tangent cones},
Ann. of Math. (2) \textbf{113} (1981), no.~2, 215--265.
\href{https://doi.org/10.2307/2006984}{doi:10.2307/2006984}.

\bibitem{CollinsLi}
T.~C. Collins and Y.~Li,
\emph{Uniqueness of some cylindrical tangent cones to special Lagrangians},
Geom. Funct. Anal. \textbf{33} (2023), 376--420.
\href{https://doi.org/10.1007/s00039-023-00634-x}{doi:10.1007/s00039-023-00634-x}.

\bibitem{DimlerGaia}
B.~Dimler and F.~Gaia,
\emph{Special Lagrangians with multiple isolated singularities},
preprint, 2026,
\href{https://arxiv.org/abs/2607.26302}{arXiv:2607.26302}.

\bibitem{EngelsteinSpolaorVelichkov}
M.~Engelstein, L.~Spolaor, and B.~Velichkov,
\emph{(Log-)epiperimetric inequality and regularity over smooth cones
for almost area-minimizing currents},
Geom. Topol. \textbf{23} (2019), no.~1, 513--540.
\href{https://doi.org/10.2140/gt.2019.23.513}{doi:10.2140/gt.2019.23.513}.

\bibitem{GilbargTrudinger}
D.~Gilbarg and N.~S.~Trudinger,
\emph{Elliptic Partial Differential Equations of Second Order},
second edition, Classics in Mathematics, Springer, Berlin, 2001.
\href{https://doi.org/10.1007/978-3-642-61798-0}{doi:10.1007/978-3-642-61798-0}.

\bibitem{HarveyLawson}
R.~Harvey and H.~B. Lawson, Jr.,
\emph{Calibrated geometries},
Acta Math. \textbf{148} (1982), 47--157.

\bibitem{Haskins}
M.~Haskins,
\emph{The geometric complexity of special Lagrangian \(T^2\)-cones},
Invent. Math. \textbf{157} (2004), no.~1, 11--70.
\href{https://doi.org/10.1007/s00222-003-0348-x}{doi:10.1007/s00222-003-0348-x}.

\bibitem{JoyceI}
D.~Joyce,
\emph{Special Lagrangian submanifolds with isolated conical singularities.
I. Regularity},
Ann. Global Anal. Geom. \textbf{25} (2004), 201--251.

\bibitem{JoyceII}
D.~Joyce,
\emph{Special Lagrangian submanifolds with isolated conical singularities.
II. Moduli spaces},
Ann. Global Anal. Geom. \textbf{25} (2004), no.~4, 301--352.
\href{https://doi.org/10.1023/B:AGAG.0000023230.21785.8d}{doi:10.1023/B:AGAG.0000023230.21785.8d}.

\bibitem{JoyceV}
D.~Joyce,
\emph{Special Lagrangian submanifolds with isolated conical singularities.
V. Survey and applications},
J. Differential Geom. \textbf{63} (2003), no.~2, 279--347.
\href{https://arxiv.org/abs/math/0303272}{arXiv:math/0303272}.

\bibitem{YangLi}
Y.~Li,
\emph{Uniqueness of tangent cones with immersed special Legendrian links},
preprint, 2026,
\href{https://arxiv.org/abs/2610.10790}{arXiv:2610.10790}.

\bibitem{Marshall}
S.~P.~Marshall,
\emph{Deformations of special Lagrangian submanifolds},
D.Phil. thesis, University of Oxford, 2002.

\bibitem{McLean}
R.~C. McLean,
\emph{Deformations of calibrated submanifolds},
Commun. Anal. Geom. \textbf{6} (1998), no.~4, 705--747.
\href{https://doi.org/10.4310/CAG.1998.v6.n4.a4}{doi:10.4310/CAG.1998.v6.n4.a4}.

\bibitem{Moriyama}
T.~Moriyama,
\emph{Deformations of special Legendrian submanifolds in Sasaki--Einstein manifolds},
Math. Z. \textbf{283} (2016), no.~3--4, 1111--1147.
\href{https://doi.org/10.1007/s00209-016-1635-8}{doi:10.1007/s00209-016-1635-8}.

\bibitem{Simon}
L.~Simon,
\emph{Asymptotics for a class of non-linear evolution equations, with
applications to geometric problems},
Ann. of Math. (2) \textbf{118} (1983), no.~3, 525--571.
\href{https://doi.org/10.2307/2006981}{doi:10.2307/2006981}.

\end{thebibliography}
\end{document}